\documentclass[11pt]{article}
\usepackage[latin1]{inputenc}
\usepackage[a4paper,width=17cm,height=24cm]{geometry}
\usepackage{amsmath,amsfonts,amssymb,indentfirst}

\usepackage{theorem}
{\theorembodyfont{\slshape}
\newtheorem{thm}{Theorem}

\newtheorem{prop}[thm]{Proposition}
\newtheorem{lemma}[thm]{Lemma}

}

\newcommand{\Proof}[1][]{{\par\smallskip\noindent{\bf Proof#1.\enspace }}}
\newcommand{\cqfd}{\hfill\rule{0.35em}{0.35em}\par\medskip\noindent}

\newcommand{\R}{\mathbb{R}}
\renewcommand{\div}{\operatorname{div}}
\newcommand{\norme}[2][]{\left\|#2\right\|_{#1}}

\usepackage{xcolor}

\title{Existence of Solutions of a\\ Non-Linear Eigenvalue Problem with a Variable Weight}
\author{Rejeb Hadiji, Fran\c{c}ois Vigneron}
\begin{document}
\maketitle

\begin{center}
Université Paris-Est,
Laboratoire d'Analyse et de Mathématiques Appliquées, UMR 8050 du CNRS\\
61, avenue du Général de Gaulle,
F-94010 Créteil -- France.
\end{center}

\begin{abstract}
We study the non-linear minimization problem on $H^1_0(\Omega)\subset L^q$ with $q=\frac{2n}{n-2}$, $\alpha>0$
and $n\geq4$~:
$$\inf_{\substack{u\in H^1_0(\Omega)\\\norme[L^q]{u}=1}}\int_\Omega a(x,u)|\nabla u|^2 - \lambda \int_{\Omega} |u|^2.$$
where $a(x,s)$ presents a global minimum $\alpha$ at $(x_0,0)$ with $x_0\in\Omega$. In order to describe the concentration of $u(x)$ around $x_0$,
one needs to calibrate the behaviour of $a(x,s)$ with respect to $s$. The model case is
$$\inf_{\substack{u\in H^1_0(\Omega)\\\norme[L^q]{u}=1}}\int_\Omega (\alpha+|x|^\beta |u|^k)|\nabla u|^2 - \lambda \int_{\Omega} |u|^2.$$
In a previous paper dedicated to the same problem with $\lambda=0$, we showed that minimizers exist only in the range $\beta<kn/q$, which corresponds
to a dominant non-linear term. On the contrary, the linear influence for
$\beta\geq kn/q$ prevented their existence.
The goal of this present paper is to show that for $0<\lambda\leq \alpha\lambda_1(\Omega)$,
$0\leq k\leq q-2$ and $\beta > kn/q + 2$, minimizers do exist.
\\[1em]
\textbf{Keywords :}
\textsl{Critical Sobolev exponent, Minimization problem, Non-linear effects.}\\
\textbf{AMS classification :} 35A01, 35A15, 35J57, 35J62.
\end{abstract}

\section{Introduction, notations and statement of the result}

\subsection{The classical non-linear problem}

The domain $\Omega$ is a smooth, bounded subset of $\R^n$ with $n\geq 4$.
Let us recall the traditional minimization problem of \cite{BN}:
\begin{equation}\label{defS}
S=\inf_{\substack{u \in H^1_0(\Omega)\\\norme[L^q]{u}=1}} \int_{\Omega} |\nabla u|^2
\end{equation}
where $q=\frac{2n}{n-2}$ is the critical exponent for the Sobolev embedding $H^1_0(\Omega)\subset L^q(\Omega)$.
For a smooth positive cut-off function $\zeta$ compactly supported and equal to 1 near the origin, a minimizing sequence for $S$ is given by~$\omega_\varepsilon/\|\omega_\varepsilon\|_{L^q}$ as $\varepsilon\to0$ with
\begin{equation}\label{omegaeps}
\omega_\varepsilon(x) = \frac{\varepsilon^{\frac{n-2}{4}}\zeta(x)}{(\varepsilon+|x|^2)^{\frac{n-2}{2}}}\cdotp
\end{equation}
According to \cite{BN}, $S$ is never achieved but one has $S=K_1/K_2$ where the constants $K_1$ and $K_2$ are the limit of the $H^1_0$ and $L^q$ norm of $\omega_\varepsilon$:
\begin{equation}
K_1 = \lim_{\varepsilon\to0}\norme[L^2]{\nabla \omega_\varepsilon}^2 \quad\text{and}\quad
K_2 = \lim_{\varepsilon\to0}\norme[L^q]{\omega_\varepsilon}^2.
\end{equation}
For example, the computation for the gradient goes as follows:
\[
\nabla\omega_\varepsilon(x)=-(n-2)
\frac{\varepsilon^\frac{n-2}{4} x\zeta(x)}{(\varepsilon+|x|^2)^\frac{n}{2}}
+\frac{\varepsilon^\frac{n-2}{4} \nabla\zeta(x)}{(\varepsilon+|x|^2)^\frac{n-2}{2}}\cdotp
\]
As $\nabla\zeta=0$ in a neighborhood of the origin,  one gets
\begin{equation}\label{gradOmega}
|\nabla\omega_\varepsilon(x)|^2
\underset{|x|\to0}{\sim} (n-2)^2 \varepsilon^\frac{n-2}{2}\frac{|x|^2 \zeta^2(x)}{{(\varepsilon+|x|^2)^n}},
\end{equation}
which integrates to a constant independent of $\varepsilon$.
In a similar way, one can deduce precise asymptotics for the various norms, that we will reuse later (see again \cite{BN}):
\begin{equation}\label{BN:gradL2andLqnorm}
\int_\Omega |\nabla \omega_\varepsilon(x)|^2dx = K_1 + O(\varepsilon^{\frac{n-2}{2}}),
\qquad
\left(\int_\Omega \omega_\varepsilon^q(x) dx\right)^{2/q} = K_2 + O(\varepsilon^{\frac{n-2}{2}}),
\end{equation}
\begin{equation}\label{BN:L2norm}
\int_\Omega \omega_\varepsilon^2(x)dx = \begin{cases}
%K_3 \sqrt{\varepsilon} + O(\varepsilon) & \text{if }n=3\\
K_3\varepsilon|\log\varepsilon|+O(\varepsilon) & \text{if }n=4\\
K_3\varepsilon + O(\varepsilon^{\frac{n-2}{2}}) & \text{if }n\geq 5.
\end{cases}
\end{equation}

\subsection{The general non-linear problem}

We are interested in the following non-linear minimization problem:
\begin{equation}\label{mainmain}
S_\lambda(a)=\inf_{\substack{u\in H^1_0(\Omega)\\\norme[L^q]{u}=1}} \left\{ \int_\Omega a(x,u)|\nabla u|^2 - \lambda \int_{\Omega} |u|^2 \right\}.
\end{equation}
where $a(x,s)=a(x,-s)$ is a smooth function, for example continuous on $\Omega\times\R$, with a continuous derivative with respect to $s$ on $\Omega\times\R_+$.
We assume that $a(x,s)$ presents a global minimum $\alpha$ at $(x_0,0)$ with $x_0\in\Omega$. One expects that minimizing sequences $u_j(t,x)$
will concentrate around~$x_0$. In order to describe the concentration of $u(x)$ around~$x_0$,
one needs to calibrate the behaviour of $a(x,s)$ with respect to $s$.

\bigskip
The model case we are focusing on is
\begin{equation}\label{generalmodel}
a(x,s)=b_1(x)+b_2(x)|s|^k
\end{equation}
where $b_1(x)\geq \alpha=b_1(x_0)$ and $b_2(x)\geq0$ with~$b_2(x_0)=0$.
Even in this simplified model, changes to the value of $k$ or of the behaviours of $b_1$ or $b_2$ around $x_0$ can lead to radically different phenomena
for the minimization problem~\eqref{mainmain}. The proper assumptions are the following:
\begin{enumerate}
\item The exponent $k$ satisfies
\begin{equation}
0\leq k < q.
\end{equation}
This range will be slightly reduced in the statement of Theorem~\ref{mainThm} and we will explain why in the next subsection.
\item $b_1$ has a global minimum $\alpha=b_1(x_0)$ at some point $x_0\in \Omega$, of order $\gamma>2$, i.e.:
\begin{equation}
\begin{cases}
b_1(x)=\alpha + O(|x-x_0|^\gamma),\qquad \gamma>2\\
b_1(x)\geq\alpha \quad\text{if}\quad  x\neq x_0\end{cases}
\end{equation}
\item$b_2$ is positive and has a unique zero in $\Omega$ at the same point $x_0$, of order $\beta\geq0$
\begin{equation}
\begin{cases}
b_2(x)=|x-x_0|^\beta + o(|x-x_0|^\beta)\\
b_2(x)>0 \quad\text{if}\quad x\neq x_0.
\end{cases}
\end{equation}
Further restrictions on the value of $\beta$ will be explained in the next subsection.
\end{enumerate}

\bigskip
For the sake of clarity, we are going to focus most of this paper on the reduced model where $b_1$ is a constant and $b_2$ is a power law.
Without restrictions, one can also assume that $x_0$ is the origin. One is thus lead to the following reduced minimization problem:
\begin{equation}\label{mainPb}
S_\lambda(\beta,k)  = \inf_{\substack{u\in H^1_0(\Omega)\\\norme[L^q]{u}=1}} E_\lambda(u)
\end{equation}
where
\begin{equation}\label{energy}
E_\lambda(u) = \int_\Omega (\alpha+|x|^\beta |u|^k)|\nabla u|^2 - \lambda \int_{\Omega} |u|^2.
\end{equation}
The case $\lambda=0$ has been extensively studied in  our previous paper, \cite{HV}.
The case $k=0$ has already been dealt with in \cite{HY}.

\bigskip
The general model case~\eqref{generalmodel} will be detailed in the last section, \S\ref{generalisation}, of this paper.
In its full generality, the qualitative properties of minimizing sequences of~\eqref{mainmain} are too varied to be described easily.

\subsubsection{Natural scaling(s) of the problem}\label{ratio}
For the sake of clarity, let us discuss the various natural possible scalings on the reduced problem~\eqref{mainPb}.
Depending on the ratio $\beta/k$, different situations occur in the blow-up scale around the point
where the weight is singular. 
More precisely, let us define $v_\varepsilon$ by $v(x) =  \varepsilon^{-n/q} v_\varepsilon(x/\varepsilon)$
and $\Omega_\varepsilon=\varepsilon^{-1}\Omega$.
\begin{itemize}
\item If $\beta < \frac{k n}{q}$ the leading term of the blow-up around $x=0$ is the non-linear one:
$$E_0 (v) \underset{\varepsilon\to0}{\sim} \varepsilon ^{-\left(\frac{kn}{q}-\beta\right)} \int_{\Omega_ \varepsilon} |y|^\beta |v_\varepsilon(y)|^k |\nabla v_\varepsilon(y)|^2 dy.$$
We showed in \cite{HV} that $E_0$ admits minimizers on $\mathbb{S}=\{u\in H^1_0(\Omega) \vert \|u\|_{L^q}=1\}$ in this case.
\item If $\beta = \frac{k n}{q}$ the first linear and non-linear terms have the same weight and
the blow-up rescaling leaves the value of $E_0(v)$ unchanged.
The corresponding infimum $S_0(\beta,k)$ does not depend on~$\Omega$ but $E_0$ admits no smooth
minimizer on $\mathbb{S}$.
\item If $\beta > \frac{kn}{q}$, the blow-up around $x=0$ gives
\begin{equation}\label{scalingoffset}
E_\lambda (v) =
\alpha \int_{\Omega_ \varepsilon} |\nabla v_\varepsilon(y)|^2 dy +
\varepsilon ^{\left(\beta-\frac{kn}{q}\right)} \int_{\Omega_ \varepsilon} |y|^\beta |v_\varepsilon(y)|^k |\nabla v_\varepsilon(y)|^2 dy
 -  \varepsilon^{2}  \lambda \int_{\Omega_\varepsilon} v_\varepsilon^2(y) dy.
 \end{equation}
\end{itemize}

When $\beta>\frac{kn}{q}$ and $\lambda=0$, we have shown in \cite{HV} that the linear behavior is dominant and
that  $E_0$ admits no minimizer on $\mathbb{S}$. One can even find a common minimizing sequences
for both the linear and the non-linear problem. A cheap way to justify this is as follows.
The problem tends to concentrate $u$ as a radial decreasing
function around the origin.
Thus, when $\beta>kn/q$, one can expect $|u(x)|^q \ll  1/|x|^{\beta q/k}$
because the right-hand side would not be locally integrable while the left-hand side is required to.
In turn, this inequality reads $|x|^\beta |u(x)|^k \ll 1$, which eliminates the non-linear contribution
in the minimizing problem $E_0$.

\medskip
When $\lambda\neq0$, the situation is quite different as both linear terms compete.
We will show in this paper that if $\beta>\frac{kn}{q}+2$, 
the problem of minimizing $E_\lambda(u)$ admits solutions. This result has the same flavor as \cite{BN} but even if the non-linear term
is not expected to be dominant, one has to deal with it rigorously.

\medskip
The gap between our existence result in \cite{HV} (namely $\beta<kn/q$) and Theorem~\ref{mainThm} below (i.e. $\beta>kn/q+2$) cannot be bridged easily.
For $\beta<\frac{kn}{q}$, we have shown in \cite{HV} that $E_0(u)$ admits minimizers. It is natural to expect that $E_\lambda(u)$ would too, as a perturbation problem.
For $\beta>\frac{kn}{q}+2$, $E_0(u)$ does not admit minimizers, but according to Theorem~\ref{mainThm} below, $E_\lambda(u)$ does for $\lambda>0$ small enough. What happens in the case of $\frac{kn}{q}\leq \beta \leq\frac{kn}{q}+2$ is not known. Let us however point out that, in the case of $\beta=2$ and $0\leq k \leq \frac{4}{n-2}$ which is a particular instance of $\frac{kn}{q}\leq \beta\leq\frac{kn}{q}+2$, the Poho\v{z}aev identity (see e.g. \cite{P} or \cite[eq. 18]{HV})
\begin{equation}\label{PZ}
\frac{1}{2} \left(\beta - \frac{kn}{q}\right) \int_\Omega |x|^{\beta} |u|^k |\nabla u|^2
+\frac{1}{2}\int_{\partial\Omega} \left|\frac{\partial u}{\partial \nu}\right|^2(x\cdot \mathbf{n})
= \lambda \int_\Omega u^2
\end{equation}
This ineqality can be restrictive for $\lambda$.
For instance, if $\Omega$ is star-shaped with respect to the origin, then $x\cdot \mathbf{n}\geq0$.
Combined with a Hardy type inequality (see for example~\cite{CKN} or \cite{HY}), it then gives us:
\begin{equation}
\frac{n^2}{8} \left(2 - \frac{kn}{q}\right) \left(\frac{k}{2}+1\right)^{-2} \int_\Omega |u|^{k+2}
\leq
\frac{1}{2} \left(2 - \frac{kn}{q}\right) \left(\frac{k}{2}+1\right)^{-2} \int_\Omega |x\cdot\nabla ( u^{\frac{k}{2}+1})|^2
\leq \lambda \int_\Omega u^2.
\end{equation}
If, for example, $k=0$, the solution $u$ does not exist if $\lambda<\frac{n^2}{4}\cdotp$

\medskip
Actually, the blow-up picture \eqref{scalingoffset} is slightly more complicated than just saying that the non-linear term disappears
because it is associated to a high power of $\varepsilon$\ldots{}
As we know, the Sobolev embedding grants $H^1_0\subset L^q\cap L^2$. But because of the identity
\[
\int_\Omega |x|^\beta |u|^k |\nabla u|^2 = \left(\frac{k}{2}+1\right)^{-2}\int_\Omega |x|^\beta (\nabla ( |u|^{\frac{k}{2}+1}))^2,
\]
the uniform boundedness of the non-linear term (as one can expect along a minimizing sequence of~\eqref{mainPb}) implies an additional
restriction on $u$, namely roughly $u^{\frac{k}{2}+1} \in L^q$ i.e. $u\in L^{\bar{q}}$ with $\bar{q}=q/(k/2+1)$. On a bounded
domain $\Omega$, this information does not seem relevant because it was already granted
by the H\"older inequality
$u\in L^q(\Omega)\subset L^{\bar{q}}(\Omega)$.
But in the blow-up process, the domain $\Omega_\varepsilon=\varepsilon^{-1}\Omega$ is ultimately rescaled to $\R^n$ and the non-linear
restriction then takes on its full significance.
For example, it could happen in~\eqref{scalingoffset} that $v_\varepsilon$ blows up to a function that does not belong
to $L^{\bar{q}}(\R^n)$. In that case, the non-linear term
\[
\varepsilon ^{\left(\beta-\frac{kn}{q}\right)} \int_{\Omega_ \varepsilon} |y|^\beta |v_\varepsilon(y)|^k |\nabla v_\varepsilon(y)|^2 dy
= C_k \varepsilon ^{\left(\beta-\frac{kn}{q}\right)} \int_{\Omega_ \varepsilon} |y|^\beta \left|\nabla (v_\varepsilon(y)^{\frac{k}{2}+1})\right|^2 dy
\]
might not be of a lower order anymore and the nature of the problem would then change completely. To avoid this potentially disastrous effect, one
needs to ensure that $L^q(\R^n)\cap L^2(\R^n)$ is embedded in $L^{\bar{q}}(\R^n)$. This restriction reads simply
$\bar{q}\geq 2$ and boils down to
\begin{equation}\label{kRestrict}
k\leq q-2.
\end{equation}
This will be part of the assumptions in Theorem~\ref{mainThm}. When on the contrary $k>q-2$, it is not clear that the non-linear term is
a lower-order term in the blow-up scaling \eqref{scalingoffset}.

\subsubsection{Admissible values for $\lambda$}\label{parRangeLambda}
In the minimization problem~\eqref{mainPb}-\eqref{energy}, the critical value for $\lambda$ is the first Dirichlet eigenvalue:
\begin{equation}\label{lambda1}
\lambda_1(\Omega)=\inf_{\substack{u\in H^1_0\\ \|u\|_{L^2}=1}} \int_\Omega |\nabla u|^2.
\end{equation}
Let us recall that $\alpha>0$ is the minimum value of the weight in our energy functional \eqref{energy}.
\begin{prop}\label{prop1}
For $0\leq \lambda\leq\alpha\lambda_1(\Omega)$ and $k>0$, one has
\begin{equation}
S_\lambda(\beta,k)\geq0.
\end{equation}
When $\lambda<\alpha\lambda_1(\Omega)$, one even has $S_\lambda(\beta,k)>0$.
\end{prop}
\Proof
Let us define:
\[
\lambda_1^{\beta,k}(\Omega) = \inf_{\substack{u\in H^1_0(\Omega)\\u\neq 0}} \frac{\int_\Omega (\alpha+|x|^\beta |u|^k)|\nabla u|^2 }{\int_\Omega u^2}\cdotp
\]
First, let us check that $\lambda_1^{\beta,k}(\Omega)=\alpha\lambda_1(\Omega)$. Indeed, let $\varphi\in H^1_0(\Omega)$
such that $-\Delta \varphi = \lambda_1(\Omega) \varphi$ with $\varphi\neq0$. Then, if one uses $\frac{1}{N}\varphi(x)$ as a test function with $N\in\mathbb{N}$, one gets, for $k>0$:
\[
0<\alpha\lambda_1(\Omega) \leq \lambda_1^{\beta,k}(\Omega) \leq \alpha\lambda_1(\Omega) + \frac{1}{N^k}\frac{\int_\Omega |x|^\beta |\varphi|^k |\nabla \varphi|^2}{\int_\Omega \varphi^2} \underset{N\to\infty}{\longrightarrow} \alpha\lambda_1(\Omega)
\]
and thus $\lambda_1^{\beta,k}(\Omega)=\alpha\lambda_1(\Omega)$.
For $u\neq0$ and $\lambda\leq\lambda_1^{\beta,k}(\Omega)$, one has therefore:
\begin{equation}\label{smalllambda}
E_\lambda(u) \geq (\lambda_1^{\beta,k}(\Omega)-\lambda)\int_\Omega u^2\geq 0,
\end{equation}
which in turn ensures that $S_\lambda(\beta,k)\geq0$.
For now, in what follows, we will not use more than this large inequality.
However, once Theorem~\ref{mainThm} is established below,
it will be asserted that $S_\lambda(\beta,k)=E_\lambda(u)$ for some non-trivial $u$ and~\eqref{smalllambda} will
then ensure that $S_\lambda(\beta,k)>0$ if $\lambda<\alpha\lambda_1(\Omega)$.
\cqfd

\bigskip
For the sake of completeness, let us briefly investigate the larger values of $\lambda$.
\begin{prop}\label{prop2}
For $\lambda\geq\alpha\lambda_1(\Omega)$ and $k>0$, one has
\begin{equation}\label{largelambda}
S_\lambda(\beta,k)\geq - (\lambda-\alpha\lambda_1(\Omega)) \cdot |\Omega|^{1-2/q}.
\end{equation}
\end{prop}
\Proof
By definition of $\lambda_1(\Omega)$, one has:
\[
E_\lambda(u) \geq \int_\Omega |x|^\beta |u|^k|\nabla u|^2 -  (\lambda-\alpha\lambda_1(\Omega)) \int_{\Omega} |u|^2.
\]
If $\|u\|_q=1$ then, according to the H\"older inequality, one has $\|u\|_2^2 \leq |\Omega|^{1-2/q}$ and \eqref{largelambda} follows immediately because $\lambda-\alpha\lambda_1(\Omega)\geq0$.
\cqfd

\paragraph{Remark}
If $\lambda>\alpha\lambda_1(\Omega)+\frac{\int_\Omega |x|^\beta |\varphi|^k
|\nabla\varphi|^2}{\|\varphi\|_q^k\|\varphi\|_2^2}$ where $\varphi$ is a non-trivial eigenfunction
of $-\Delta \varphi = \lambda_1(\Omega) \varphi$ then,  using 
$u=\varphi/\|\varphi\|_q$ as a test function, one gets $E_\lambda(u)<0$ and thus $S_\lambda(\beta,k)< 0$.

\subsection{Statement of the main result}

In this paper, our main result is the following.

\begin{thm}\label{mainThm}
Let us assume that $n\geq4$. One sets $q=\frac{2n}{n-2}>2$ and assumes that
 \begin{equation}
 0< \lambda\leq \alpha\lambda_1(\Omega), \qquad
0\leq k\leq q-2 \qquad\text{and}\qquad \beta > \frac{kn}{q}+2.
\end{equation}
Then there exists $u\in H^1_0(\Omega)$ with $\|u\|_{L^q}=1$ such that $E_\lambda(u)=S_\lambda(\beta,k)$.
\end{thm}
In section \S\ref{generalisation}, one will prove a similar result about the general model~\eqref{generalmodel},
which will be stated as Theorem~\ref{mainThm2}.

\medskip\noindent
Over the course of the proof, one also reaps the following convergence result.
\begin{prop}\label{corollary}
Under the same assumptions,
for any minimizing sequence $u_j \in H^1_0(\Omega)$ i.e. such that
\[
\|u_j\|_{L^q}=1 \quad\text{and}\qquad E_\lambda(u_j)=S_\lambda(\beta,k)+o(1)
\]
that converges weakly to some $u\in H^1_0(\Omega)$, then $\|u\|_{L^q}=1$ and
the convergence also holds in the strong topology of $H^1(\Omega)$.
Moreover, $u$ is a minimizer and solves the Euler-Lagrange equation:
\begin{equation}\label{eqEL}
\begin{cases}
-\div\left(\left(\alpha+|x|^\beta |u|^k\right)\nabla u \right)+ \frac{k}{2} |x|^\beta |u|^{k-2}u |\nabla u|^2 =\lambda u + \Theta |u|^{q-2}u\\
u_{\vert \partial\Omega}=0
\end{cases}
\end{equation}
for some $\Theta>0$.
\end{prop}
Let us point out that the uniqueness of the limit is not known and constitutes a wide-open problem.
For example, one knows at least that uniqueness does not hold in some cases involving weights that concentrate
on multiple origins \cite{HMPY}.

\medskip
Let us also note that if $u_j$ is a minimizing sequence, then it is standard to check that
$|u_j|$ is also a minimizing sequence. Thanks to Proposition~\ref{corollary} its limit $|u|$ is a positive
nontrivial solution
of \eqref{eqEL} with an $L^q$-norm equal to 1.

\bigskip\noindent%
Let us now comment briefly upon the assumptions of the main Theorem.
\begin{itemize}
\item
The restriction $0< \lambda\leq \alpha\lambda_1(\Omega)$ is natural in regard to Proposition~\ref{prop1}
and ensures $S_\lambda(\beta,k)\geq0$.
For partial results in the case $\lambda>\alpha\lambda_1(\Omega)$, see the concluding remark
of sections \S\ref{parRangeLambda} and \S\ref{proofMainThm}.
\item
The restriction $\beta>\frac{kn}{q}+2$ comes from the competition between the different scalings of
the terms that appear in the expression of $E_\lambda(u)$. In particular, the $+2$ offset reflects the
scaling of $\int_\Omega u^2$, as pointed out in the discussion of the previous section, \S\ref{ratio} about~\eqref{scalingoffset}.
It is a crucial assumption that provides $S_\lambda(\beta,k)<\alpha S$, which will be shown in Lemma~\ref{aprioriestimate} below.
\item
The restriction $k\leq q-2$ has also been discussed in section \S\ref{ratio}. It is necessary to ensure that the non-linear term
stays of lower order throughout the minimization process, especially  at the finest scales around the singularity $x=0$.
Let us also observe that for the critical exponent $k=q-2$, both terms $-\div( |x|^\beta |u|^k \nabla u)$ and $|u|^{q-2}u$
of the Euler-Lagrange equation have the same weight for the amplitude scaling transform $u\mapsto A u$ with $A>0$. 
\item
The case $\lambda=0$ has been extensively studied in \cite{HV} and the behavior for $\beta>\frac{kn}{q}$
is then quite different from that in Theorem~\ref{mainThm} because the non-existence of minimizers was established.
\item
It is also natural to exclude the case of $\lambda<0$.
Indeed, let us consider the classical  Poho\v{z}aev identity~\eqref{PZ} when $\Omega$ is star-shaped.
For $\beta\geq  \frac{kn}{q}$, the left-hand side of \eqref{PZ} is positive but the right-hand side is negative if $\lambda< 0$.
There are therefore no minimizers when $\lambda<0$, at least for star-shape domains.
\item
When $k=0$, the non-linear nature of the problem changes. The critical value is then $\beta=2$ and is excluded. In that case, Theorem~\ref{mainThm} contains the results obtained in \cite{HY}.
\end{itemize}

\bigskip
Our proof of Theorem~\ref{mainThm} follows the general principles of the method of concentration~\cite{BN}. The core of our argument is a standard calculus of variation
around the weak limit  $u$ of a minimizing sequence, either in the direction of $u$ itself or in the direction of $\omega_\varepsilon$.

Thirty years later, this method can be seen
as the ``pedestrian way'' that goes along the ``highway'' of the more general method of concentration-compactness using a profile decomposition \cite{Struwe}, \cite{Lions}.
Let us give a rough sketch of this last method and point out the specificity of our problem that led us to choosing this rather~``historical'' approach.

Given a minimizing sequence $u_k\in H^1_0(\Omega)\cap\{ v \,\vert\, \|v\|_{L^q}=1\}$ of some functional $F$ (in our case $F(u)=\int a(x,u)|\nabla u|^2-\lambda\int u^2$),
the boundness of the sequence in $H^1_0(\Omega)$ allows us to assume that, up to some sub-sequence, the sequence $u_k$ converges weakly to some limit function $u\in H^1_0(\Omega)\cap\{ v \,\vert\, \|v\|_{L^q}\leq1\}$.
The question of showing that $u$ is a minimizer is roughly equivalent to showing that the convergence of this sequence holds in the strong topology of $L^q(\Omega)$ and
this last statement is, at least, clearly sufficient. In case $u_k\not\to u$ in $L^q(\Omega)$, the concentration-compactness principle \cite{Sol}, \cite{gerard}, \cite{jaffard}, \cite{AT},  \cite{ST} provides orthogonal profiles
\[
u_k= u+\sum R_j(w_j) + \zeta_k \qquad\text{with}\qquad  \zeta_k \overset{L^q}{\longrightarrow} 0,
\]
where $R_j$ are translation-scaling operators at either different scales and/or locations and $w_j$ are calibrated profiles. 
Applying an appropriate generalized Brezis-Lieb lemma \cite{bl}, \cite{AT2}, \cite{ST}, \cite{EM} would then provide
\[
F(u_k)= F(u)+\sum F(R_j(w_j)) + o(1).
\]
In most semilinear problems, the functional $F$ is invariant by translation and scaling, so the action of $R_j$ would commute easily with the one of $F$ and one would then get
\[
F(R_j(w_j)) \geq S_0
\]
where $S_0$ is the similar minimization problem on the whole space $\R^n$. From the a-priori comparison $0\leq \inf F < S_0$ one could then deduce that the profiles just don't exist
i.e. that $u_k\to u$ strongly in $L^q$.  This comparison between $\inf F$ and the corresponding whole space problem is specific to each functional at hand and will be done here in Lemma~\ref{aprioriestimate}.

The specificity of our problem is that 
the weight $a(x,u)$ depends on the space variable in a non-trivial way. The whole point of our paper is to show a specific behaviour around a global minimum of $a$, under
some structure assumptions. We are not convinced that one of the aforementioned generalized Brezis-Lieb lemma would indeed simplify the problem to the point of triviality.

On the contrary, a precise computation of the concentration of $u_k$ around the minimum of the weight $a(x,u)$ in the spirit of \cite{BN} might, at first sight, appear outdated.
But, in our present case, it is more enlightening and it will thus be our course of action.

\subsection{Motivations and related questions}

The study of the general problem~\eqref{mainmain} is related to the associated Euler-Lagrange PDE which reads formally:
\begin{equation}\label{formalEL}
\begin{cases}
-\div\left( a(x,u) \nabla u\right) +  \frac{1}{2} \partial_s a(x,u) |\nabla u|^2 =\lambda u + \Theta
|u|^{q-2}u\\
u_{\vert \partial\Omega}=0.
\end{cases}
\end{equation}
As the weight $a(x,s)$ depends in a non-trivial way in $x$ and $u$, this equation is of quasi-linear type. It only boils down to a semi-linear problem when $\partial_s a=0$
(i.e. when $k=0$ in the reduced model~\eqref{energy}), which is only one very particular case among the general assumptions of Theorem~\ref{mainThm}.

\bigskip
The most recents developments on quasi-linear elliptic equations seem to be focused on $p$-laplace operators i.e. a leading operator of the form
$-\div(|\nabla u|^{p-2}\nabla u)$.
Recent papers that include a critical non-linearity (see \cite{DH}, \cite{dVS}, \cite{FC}, \cite{DM})
seem to either focus on singular weights that can be controlled by a Hardy-type inequality~\cite{CKN}:
\[
\int_\Omega \frac{|u(x)|^{q^\ast}}{|x|^s} dx \leq C_s \|u\|_{H^1_0(\Omega)}^{q^\ast}
\]
or focus on a term $|u|^{p^\ast-2}u$ where $p^\ast=pn/(n-p)$ is the critical exponent associated to the $p$-Laplacian.
In the case of \cite{DM}, the framework is that of the Heisenberg group.

\bigskip
The quasilinear operator $-\div(a(x,u)\nabla u)$ has been studied in various papers, however not in the critical case. 
The paper \cite{ABPP} deals with a general quasilinear elliptic equation of the form
\[
-\div(a(x,u)\nabla u) + g(x,u,\nabla u) = \lambda h(x) u + f
\]
with a quadratic growth of $g(x,u,\nabla u)$ with respect to $\nabla u$, and very general weights $a$, $g$ and $h$. However, this equation does not contain a critical non-linearity
like \eqref{formalEL}. 
In \cite{B}, the form of the quasi-linear operator gets closer to ours, but corresponds only to the case $\beta=0$:
\[
-\div((a(x)+|u|^q)\nabla u) + b(x) u|u|^{p-1}|\nabla u|^2 = f \in L^m.
\]
This problem is truly non-linear because $p>2$ but it remains sub-critical because $m<\frac{2n}{n+2}$.
In \cite{BP}, the critical points of the functional
\[
I(u)=\frac{1}{2}\int_\Omega a(x,u) |\nabla u|^2 - \frac{1}{p}\int_\Omega u_+^p
\]
are studied
for $p<2n/(n-2)$ and a general function $a$. Again, it is a sub-critical non-linearity.

\bigskip
The study of the minimization problem \eqref{mainmain} in general and of the particular form of the energy functionnal~\eqref{energy}
is motivated by deeper questions in geometry. In \cite{DN}, the minimization problem in the case $\lambda=0$ is studied on compact manifolds, along with various variants.

More generally, the Dirichlet energy associated to a map $u:(M,g)\to(N,h)$ between Riemanian manifolds takes the form
\[
\mathcal{E}(u)=\int_M g^{ij}(x) h_{kl}(u) \partial_i u^k \partial_j u^l (\det |g|)^{n/2} dx^1\ldots dx^n.
\]
Our scalar problem with $\mathcal{E}(u)=\displaystyle \int_{\R^n} a(x,|u|) |\nabla u|^2$ and $a(x,s)=\alpha+|x|^\beta s^k$ can be seen as an
entry-level model for the more general case, but with a singular metric $h$ and a non-linear term
\[
\int_\Omega |x|^\beta |u|^k |\nabla u|^2 = \left(\frac{k}{2}+1\right)^{-2}\int_\Omega |x|^\beta (\nabla ( |u|^{\frac{k}{2}+1}))^2
\]
We refer to \cite{jost} for a general survey of those questions.

\medskip
Other geometric motivations, in particular in relation to the Yamabe problem, can be found for example in \cite{A}, \cite{CNV} \cite{LP},  \cite{Lions}.
Note that the shape of the domain can have a strong influence on the type of results one can expect. See for example the
seminal work of J.M.~Coron \cite{c}, or \cite{h}.

\bigskip
Let us finally point out that the dimension $n=3$ could also be interesting for this problem, but is not yet fully understood. See \cite{BN}, \cite{CHL}.

\subsection{Structure of the paper}

The paper is structured as follows.
The next section, \S\ref{par:apriori}, proves the a-priori estimate $S_\lambda(\beta,k)<\alpha S$.
Theorem~\ref{mainThm}, which is the main result of this paper,
is proved by a contradiction argument that spans the whole of section \S\ref{par:contradict}. In section \S\ref{par:notzero} one shows that the weak limit $u$ of a minimizing
sequence is not zero. In sections \S\ref{parcv1} and \S\ref{parcv2}, we carry out a calculus of variation around
$u$, respectively in the direction of $\omega_\varepsilon$ defined by \eqref{omegaeps} and then along
$u$ itself. The sections \S\ref{proofMainThm} and \S\ref{proofMainThm2} then put the proof together.

In the final section \S\ref{generalisation}, one explains how our proof could be adapted to deal with a more general minimization problem, which is stated in
Theorem~\ref{mainThm2}.

\section{A priori estimate on $S_\lambda(\beta,k)$}\label{par:apriori}

Taking $\|\omega_\varepsilon\|_{L^q}^{-1}\omega_\varepsilon$ as a test function in $E_\lambda(u)$ provides a natural upper bound for $S_\lambda(\beta,k)$:

$$S_\lambda(\beta,k) \leq E_\lambda\left(\|\omega_\varepsilon\|_{L^q}^{-1}\omega_\varepsilon\right).$$
We will show the following:
\begin{lemma}\label{aprioriestimate}
For $n\geq4$, $\lambda>0$, $k\geq0$ and $\beta > \frac{kn}{q}+2$, one has:
\begin{equation}
S_\lambda(\beta,k) < \alpha S
\end{equation}
where $\alpha$ is the minimum value of the weight in the energy functional \eqref{energy}.
\end{lemma}
This lemma is responsible for the main limitation on $\beta$ in Theorem~\ref{mainThm}.

\Proof
One needs to compute precisely the asymptotic expansion with respect to $\varepsilon$ of:
\[
E_\lambda\left(\frac{\omega_\varepsilon}{\|\omega_\varepsilon\|_q}\right)
= E_0\left(\frac{\omega_\varepsilon}{\|\omega_\varepsilon\|_q}\right)
-\frac{\lambda}{\|\omega_\varepsilon\|_q^2}\int_\Omega \omega_\varepsilon^2
\]
Each term in the expression $ E_0\left(\frac{\omega_\varepsilon}{\|\omega_\varepsilon\|_q}\right)$
has been studied in our previous paper, \cite{HV}:
\[
E_0\left(\frac{\omega_\varepsilon}{\|\omega_\varepsilon\|_q}\right) =
\frac{\alpha}{\|\omega_\varepsilon\|_q^2}\int_\Omega |\nabla \omega_\varepsilon|^2
+ \frac{1}{\|\omega_\varepsilon\|_q^{k+2}}\int_\Omega |x|^\beta |\omega_\varepsilon|^k |\nabla \omega_\varepsilon|^2
\]
with, according to \eqref{BN:gradL2andLqnorm}:
\[
\frac{1}{\|\omega_\varepsilon\|_q^2}\int_\Omega |\nabla \omega_\varepsilon|^2 =
\frac{ K_1  + O(\varepsilon^{\frac{n-2}{2}})}{ K_2 + o(\varepsilon^{\frac{n-2}{2}})}=S+O(\varepsilon^{\frac{n-2}{2}})
\]
and (see \cite[Proposition 5]{HV}):
\[
\frac{1}{\|\omega_\varepsilon\|_q^{k+2}}\int_\Omega |x|^\beta |\omega_\varepsilon|^k |\nabla \omega_\varepsilon|^2 = %\frac{1}{K_2^{(k+2)/2} + o(\varepsilon^{\frac{n-2}{2}})}
\left\{\begin{array}{llll}
O\left(\varepsilon^{\frac{2\beta-k(n-2)}{4}}\right)&\textrm{\,if\enspace
$\frac{kn}{q}<\beta<(k+1)(n-2)$}\\[\medskipamount]
O\left(\varepsilon^{\frac{(k+2)(n-2)}{4}}\:|\log\varepsilon|\right)&\textrm{\,if
$\beta =(k+1)(n-2)$}\\[\medskipamount]
O\left(\varepsilon^{\frac{(k+2)(n-2)}{4}}\right)&\textrm{\,if $\beta >
(k+1)(n-2).$}
\end{array}\right.
\]
%For $n\geq4$ and $\beta>\frac{kn}{q}+2$, one has necessarily $\beta>(k+1)(n-2)$. Indeed
%\[
%(k+1)(n-2)\geq\frac{kn}{q}+2 \quad\Longleftrightarrow\quad k\geq -\frac{2(n-4)}{n-2}
%\]
Note that $(k+1)(n-2)=\frac{kn}{q}+\frac{(k+2)n}{q}$.
Thanks to \eqref{BN:L2norm}, the additional term satisfies:
\[
-\frac{\lambda}{\|\omega_\varepsilon\|_q^2}\int_\Omega \omega_\varepsilon^2
= -\frac{\lambda K_3}{K_2}%\left(1+o(\varepsilon^{\frac{n-2}{2}})\right)\times 
 \begin{cases}
%\sqrt{\varepsilon} + O(\varepsilon) & \text{if }n=3\\
\varepsilon|\log\varepsilon|+O(\varepsilon) & \text{if }n=4\\
\varepsilon + O(\varepsilon^{\frac{n-2}{2}}) & \text{if }n\geq 5.
\end{cases}
\]
This term is clearly the dominant remainder and dictates the sign if $k>0$ and $\beta>\frac{kn}{q}+2$. Indeed, one has for $n\geq4$:
\begin{equation}\label{equivTech}
\min\left\{\frac{2\beta - k(n-2)}{4} ; \frac{(k+2)(n-2)}{4}\right\}>1
\quad\Longleftrightarrow\quad
\begin{cases}
\beta>\frac{kn}{q}+2 & \text{if } \beta < \frac{kn}{q}+\frac{(k+2)n}{q}\\
k>\frac{-2(n-4)}{n-2} & \text{if } \beta \geq \frac{kn}{q}+\frac{(k+2)n}{q}
\end{cases}
\end{equation}
and thus
%(note that $\log \varepsilon$ is only necessary in the second term when $\beta=(k+1)(n-2)$, but
%here, one can keep it for any $\beta$ anyway):
\[
E_\lambda\left(\frac{\omega_\varepsilon}{\|\omega_\varepsilon\|_q}\right)=
\alpha S+O(\varepsilon^{\frac{n-2}{2}})+
O\left(\varepsilon^{\min\left\{\frac{2\beta - k(n-2)}{4} ; \frac{(k+2)(n-2)}{4}\right\}}|\log \varepsilon|\right)
-\frac{\lambda K_3}{K_2}
\begin{cases}
\varepsilon|\log\varepsilon|+O(\varepsilon) & \text{if }n=4\\
\varepsilon + O(\varepsilon^{\frac{n-2}{2}}) & \text{if }n\geq 5.
\end{cases}
\]
If $k=0$, one needs to distinguish between dimensions. If $n\geq 5$, the $-\lambda$ term is dominant because
$\frac{-2(n-4)}{n-2}<0$. If $n=4$, and $\beta>2$, the equivalence \eqref{equivTech} cannot help anymore but
the previous estimates directly give:
\[
E_\lambda\left(\frac{\omega_\varepsilon}{\|\omega_\varepsilon\|_q}\right) = \alpha S+O(\varepsilon)+O(\varepsilon)-\frac{\lambda K_3}{K_2} \varepsilon |\log \varepsilon|
\]
so $-\lambda$, again, dictates the sign.
\cqfd

\paragraph{Remark} If $\beta=\frac{kn}{q}+2$, the non-linear term has exactly the same weight as the $-\lambda$ term  so the results holds if $\lambda$ is large enough, namely
\[
C -\lambda \frac{K_3}{K_2}<0  \qquad i.e\qquad \frac{C K_2}{K_3}<\lambda.
\]
Sadly, the comparison between this critical value and $\alpha \lambda_1(\Omega)$ is not known. Note that the previous result would hold for a nontrivial range of $\lambda$ if $\alpha$ is large enough.

\section{Existence of minimizers}\label{par:contradict}
Let us take a minimizing sequence $(u_j)_{j\in\mathbb{N}}$ for  $S_\lambda(\beta,k)$.
It is a bounded sequence in $H^1_0(\Omega)$ and one can therefore consider $u\in H^1_0(\Omega)$ a weak limit of a suitable subsequence with
\begin{equation}\label{deft}
t=\|u\|_{L^q} \in [0,1]
\end{equation}
and $u_j\to u$ strongly in $L^p(\Omega)$ for any $p\in [2,q)$.
Let us assume, by contradiction, that $t<1$. \label{beginContradiction}

\subsection{The weak limit is not identically zero}\label{par:notzero}
%As the minimization problem $S$ is never achieved,
By definition~\eqref{defS} of $S$, one has
$$\alpha S \leq \alpha \int_\Omega |\nabla u_j|^2 = S_\lambda(\beta,k) +\lambda \int_\Omega u_j^2 - \int_\Omega |x|^\beta |u_j|^k |\nabla u_j|^2 + o(1).$$
Discarding the negative term and using $u_j\to u$ strongly in $L^2$ gives
$$\alpha S \leq S_\lambda(\beta,k) +\lambda \int_\Omega u^2 + o(1)$$
and thus, according to Lemma~\ref{aprioriestimate}:
$$\lambda \int_\Omega u^2 \geq \alpha S-S_\lambda(\beta,k)>0,$$
which ensures that $u$ is not identically zero, and in particular that $t\neq0$ (with $t$ defined by~\eqref{deft}).

\subsection{First calculus of variations around the weak limit}\label{parcv1}

In this section, let us explore $H^1_0(\Omega)$ around the weak limit $u$ using $u+\theta \omega_\varepsilon$ as a test function and with~$\theta$ chosen to satisfy the $L^q$-norm constraint.
\begin{lemma}\label{lemma4}
For any $v\in H^1_0(\Omega)$ such that $\|v\|_q\leq1$  and $\int_\Omega |x|^\beta |v|^k |\nabla v|^2<\infty$,
one has:
\begin{equation}\label{identity0}
S_\lambda(\beta,k) \leq E_\lambda(v) +\alpha  S \{1-\|v\|_q^q\}^{2/q}.
\end{equation}
For the weak limit $u$ of the above minimizing sequence $(u_j)_{j\in\mathbb{N}}$, one gets an equality in~\eqref{identity0}:
\begin{equation}\label{identity1}
S_\lambda(\beta,k) = E_\lambda(u) + \alpha S(1-t^q)^{2/q}
\end{equation}
with $t=\|u\|_q$ defined by~\eqref{deft}.
\end{lemma}

\Proof[ of \eqref{identity0}] 
If $\|v\|_q=1$, then $v$ is an admissible test function and \eqref{identity0} holds by definition.
Let us now assume that $\|v\|_q<1$. For each $\varepsilon>0$, the intermediary values theorem
ensures the existence of $c_\varepsilon>0$ such that $\|v+c_\varepsilon \omega_\varepsilon\|_q=1$.
One has therefore:
\[
S_\lambda(\beta,k) \leq E_\lambda(v+c_\varepsilon \omega_\varepsilon).
\]
The Brezis-Lieb lemma \cite{bl} allows one to compute $c_\varepsilon$:
\[
1=\|v\|_q^q + c_\varepsilon^q \|\omega_\epsilon\|_q^q + o(1)
\]
thus
\begin{equation}\label{cepsilon}
c_\epsilon^2 = \frac{S(1-\|v\|_q^q)^{2/q}}{K_1}+o(1).
\end{equation}
Next, one computes $E_\lambda(v+c_\varepsilon \omega_\varepsilon)$~:
\[
S_\lambda(\beta,k)\leq \int_\Omega (\alpha +|x|^\beta|v+c_\varepsilon \omega_\varepsilon|^k)|\nabla(v+c_\varepsilon \omega_\varepsilon)|^2 -\lambda \int_\Omega |v+c_\varepsilon \omega_\varepsilon|^2.
\]
Thus one has
\[
S_\lambda(\beta,k)\leq E_\lambda(v) + \alpha c_\varepsilon^2 \int_\Omega |\nabla \omega_\epsilon|^2+R_\varepsilon
\]
with $R_\varepsilon = R_\varepsilon^\text{sub}+R_\varepsilon^\text{crit}$ and
\begin{align*}
R_\varepsilon^\text{sub} &= 2 \alpha c_\varepsilon \int \nabla v\cdot \nabla\omega_\varepsilon
-2 c_\varepsilon\lambda \int_\Omega v \omega_\varepsilon 
-\lambda c_\epsilon^2\int_\Omega \omega_\epsilon^2\\
R_\varepsilon^\text{crit} &=
\int_\Omega |x|^\beta|v+c_\varepsilon \omega_\varepsilon|^k |\nabla(v+c_\varepsilon \omega_\varepsilon)|^2
-\int_\Omega |x|^\beta |v|^k |\nabla v|^2
\end{align*}
and the whole point is to show that $R_\varepsilon=o(1)$ as $\epsilon\to0$.

\medskip
For the first term of $R_\varepsilon^\text{sub}$,
one uses simply that $\omega_\varepsilon\rightharpoonup0$ weakly in $H^1_0(\Omega)$
and strongly in $L^2(\Omega)$. For the next two terms of $R_\varepsilon^\text{sub}$, one combines \eqref{BN:L2norm} and \eqref{cepsilon} to get
$O(\|\omega_\varepsilon\|_2)+O(\|\omega_\varepsilon\|_2^2)$ thus indeed $R_\varepsilon^\text{sub}$
converge to zero with $\varepsilon$.

\medskip
All that remains is to study:
\begin{align*}
R_\varepsilon^\text{crit} &=
\int_\Omega |x|^\beta \left( |v+c_\varepsilon \omega_\varepsilon|^k  -|v|^k \right) |\nabla v|^2\\
&\qquad
+ c_\varepsilon^2 \int_\Omega |x|^\beta|v+c_\varepsilon \omega_\varepsilon|^k |\nabla \omega_\varepsilon|^2
+2 c_\varepsilon \int_\Omega |x|^\beta|v+c_\varepsilon \omega_\varepsilon|^k \nabla v \cdot \nabla \omega_\varepsilon.
\end{align*}
The key is the following identity on $\R^2$:
\begin{equation}\label{splitId}
\left||x+y|^k-|x|^k\right| \leq \begin{cases}
|y|^k & \text{if } 0\leq k \leq 1,\\
|y|^k+C_k(|x|^{k-1}|y|+|x||y|^{k-1}) & \text{if } k>1,
\end{cases}
\end{equation}
which follows respectively from $||1+t|^k-1|\leq |t|^k$ if $k\in[0,1]$ and
$\left| \frac{|1+t|^k - 1 - |t|^k}{(1+|t|^{k-2})t} \right|\leq C_k$ if $k\geq0$, applied for $t=y/x\in\R$.
Let us compute those two limits first:
\begin{gather}
\int_\Omega |x|^\beta \omega_\varepsilon^k |\nabla v|^2 = o(1), \label{idA}
\\
\int_\Omega |x|^\beta (|v|+c_\varepsilon \omega_\varepsilon)^k |\nabla \omega_\varepsilon|^2=o(1). \label{idB}
\end{gather}
Subsequently, one will also check that:
\begin{gather}
\int_\Omega |x|^\beta|v+c_\varepsilon \omega_\varepsilon|^k \nabla v \cdot \nabla \omega_\varepsilon=o(1), \label{idC}\\
\int_\Omega |x|^\beta \omega_\varepsilon^{k-1}  |v| |\nabla v|^2 = o(1) \qquad \text{if $k>1$}, \label{idD}\\
 \int_\Omega |x|^\beta v^{k-1} \omega_\varepsilon |\nabla v|^2 = o(1)\qquad \text{if $k>1$}. \label{idE}
\end{gather}
Once this verification is complete, one can ascertain that $R_\varepsilon^\text{crit}$ converges
to zero with $\varepsilon$ and the proof of~\eqref{identity0} will therefore be complete.

\Proof[ of \eqref{idA}]
One uses an $L^\infty\times L^1$ estimate:
\[
\int_\Omega |x|^\beta \omega_\varepsilon^k |\nabla v|^2
\leq C\left(\int_\Omega |\nabla v|^2\right) \times \left(\sup_{r\leq\delta}  \frac{r^{\beta}  \varepsilon^{\frac{k(n-2)}{4}}}{(\varepsilon+r^2)^{\frac{k(n-2)}{2}}} \right)
\]
with $\delta =\sup_{x\in\Omega} |x|$. The right-hand side is maximal around $r\sim\sqrt{\varepsilon}$
and its maximal value is of order~$\varepsilon^{\frac{\beta}{2}-\frac{k(n-2)}{4}}$, which tends to zero
provided
\[
\beta > \frac{kn}{q}\cdotp
\]
\Proof[ of \eqref{idB}]
One uses H\"older's $L^{q/k}\times L^{q/(q-k)}$ inequality with 
$|v+c_\varepsilon \omega_\varepsilon|^k \in L^{q/k}$:
\[
\int_\Omega |x|^\beta|v+c_\varepsilon \omega_\varepsilon|^k |\nabla \omega_\varepsilon|^2
\leq
\left( \|v\|_{L^q} + \|\omega_\varepsilon\|_{L^q}\right)^k \times
\left( \int_\Omega |x|^\frac{\beta q}{q-k}|\nabla\omega_\varepsilon|^\frac{2q}{q-k} dx\right)^{1-k/q}.
\]
The first factor is bounded.
The precise computation of the gradient \eqref{gradOmega} provides the necessary decay:
\[
\left( \int_\Omega |x|^\frac{\beta q}{q-k}|\nabla\omega_\varepsilon|^\frac{2q}{q-k} dx\right)^{1-k/q}
\leq C \varepsilon^\frac{n-2}{2} \left(
\varepsilon^{-\frac{q(n-2)}{q-k}}
\int_0^\varepsilon r^{\frac{\beta q}{q-k}+n-1} dr
+\int_\varepsilon^\delta
\frac{r^{\frac{(\beta+2)q}{q-k}+n-1}}{(\varepsilon+r^2)^\frac{qn}{q-k}} dr\right)^{1-k/q}.
\]
The first term (small scale) is due to the cut-off function but is harmless
because it ultimately boils down to $\varepsilon^{\beta-\frac{kn}{q}+n\left(1-\frac{1}{q}\right)}$.
The second one is dealt with using a blow-up rescaling ($r=\sqrt{\varepsilon}\rho$):
\[
\left( \int_\Omega |x|^\frac{\beta q}{q-k}|\nabla\omega_\varepsilon|^\frac{2q}{q-k} dx\right)^{1-k/q}
\leq C \left\{
\varepsilon^{\frac{n+2}{2}+\beta-\frac{kn}{q}}
+\varepsilon^{\frac{n-2}{2}+\frac{1}{2}\left(\beta-\frac{kn}{q}\right)}
\left( \int_0^{\delta/\sqrt{\varepsilon}} 
\frac{\rho^{\frac{(\beta+2)q}{q-k}+n-1}}{(1+\rho^2)^\frac{qn}{q-k}} d\rho
\right)^{1-k/q}
\right\}
\]
There are three cases:
\begin{itemize}
\item If $\frac{kn}{q}<\beta<\frac{kn}{q}+n-2$, the last integral can be extended over $\R_+$ and gives a harmless constant factor.
\item If $\beta=\frac{kn}{q}+n-2$, the last integral is of order $\log(\delta/\sqrt{\varepsilon})$ but
as $n\geq3$, it does not prevent the whole term from tending to zero with $\varepsilon$.
\item If $\beta>\frac{kn}{q}+n-2$, the last integral boils down to
\[
\left( \int_0^{\delta/\sqrt{\varepsilon}} 
\frac{\rho^{\frac{(\beta+2)q}{q-k}+n-1}}{(1+\rho^2)^\frac{qn}{q-k}} d\rho
\right)^{1-k/q} \leq C \varepsilon^{\beta-\frac{kn}{q}+2} = O(\varepsilon^n).
\]
\end{itemize}
In all cases, assertion \eqref{idB} holds true.

\Proof[ of \eqref{idC}, assuming $0\leq k \leq 1$]
Thanks to \eqref{splitId}, it is sufficient to control the following integrals:
\[
\int_\Omega |x|^\beta|v|^k |\nabla v|| \nabla \omega_\varepsilon|
\leq \left(\int_\Omega |x|^\beta |v|^k |\nabla v|^2 \right)^{1/2} \times
\left( \int_\Omega |x|^\beta |v|^k |\nabla \omega_\varepsilon|^2  \right)^{1/2}
\]
whose factors are respectively bounded by assumption and controlled by \eqref{idB}, and
\[
\int_\Omega |x|^\beta \omega_\varepsilon^k \nabla v \cdot \nabla \omega_\varepsilon
\leq \left(\int_\Omega |x|^\beta \omega_\varepsilon^k |\nabla v|^2  \right)^{1/2} \times
\left( \int_\Omega |x|^\beta \omega_\varepsilon^k |\nabla \omega_\varepsilon|^2 \right)^{1/2} 
\]
whose factors are respectively controlled by \eqref{idA} and \eqref{idB}.
This proves \eqref{idC} when $k\in[0,1]$.

\Proof[ of \eqref{idD}, assuming $k>1$]
Let us set $\vartheta=1/k\in]0,1[$.
One uses H\"older's $L^{k}\times L^{k/(k-1)}$ inequality with 
$(|x|^\beta|\nabla v|^2)^\vartheta |v| \in L^k$ and $(|x|^\beta|\nabla v|^2)^{1-\vartheta} \omega_\varepsilon^{k-1} \in L^{k/k-1}$:
\[
\int_\Omega |x|^\beta \omega_\varepsilon^{k-1}  |v| |\nabla v|^2
\leq \left(\int_\Omega |x|^\beta |v|^k |\nabla v|^2 \right)^{1/k} \times \left(
\int_\Omega |x|^\beta \omega_\varepsilon^k |\nabla v|^2\right)^{1-1/k}.
\]
The first integral is bounded by assumption and the second one is controlled by \eqref{idA}.

\Proof[ of \eqref{idE}, assuming $k>1$]
Again, one uses H\"older's $L^{k}\times L^{k/(k-1)}$ inequality but this time with 
${(|x|^\beta|\nabla v|^2)^\vartheta \omega_\varepsilon \in L^k}$ and $(|x|^\beta|\nabla v|^2)^{1-\vartheta} v^{k-1} \in L^{k/k-1}$:
\[
\int_\Omega |x|^\beta v^{k-1} \omega_\varepsilon |\nabla v|^2
\leq  \left(\int_\Omega |x|^\beta \omega_\varepsilon^k |\nabla v|^2 \right)^{1/k} \times \left(
\int_\Omega |x|^\beta v^k |\nabla v|^2\right)^{1-1/k}.
\]
The first integral is controlled by \eqref{idA} and the second one is bounded by assumption.

\Proof[ of \eqref{idC}, assuming $k>1$]
Again, one uses \eqref{splitId} to split the integral:
\begin{align*}
\int_\Omega |x|^\beta|v+c_\varepsilon \omega_\varepsilon|^k |\nabla v||\nabla \omega_\varepsilon| &\leq
\int_\Omega |x|^\beta |v|^k |\nabla v||\nabla \omega_\varepsilon| + c_\varepsilon^k
\int_\Omega |x|^\beta \omega_\varepsilon^k |\nabla v||\nabla \omega_\varepsilon|\\
&\qquad + C_{k,\varepsilon} \left\{
\int_\Omega |x|^\beta |v|^{k-1}\omega_\varepsilon |\nabla v||\nabla \omega_\varepsilon|
+ \int_\Omega |x|^\beta v \omega_\varepsilon^{k-1} |\nabla v||\nabla \omega_\varepsilon|
\right\}.
\end{align*}
Each integral can now be controlled using Cauchy-Schwarz and the previous inequalities, namely:
\begin{align*}
\int_\Omega |x|^\beta|v+c_\varepsilon \omega_\varepsilon|^k |\nabla v||\nabla \omega_\varepsilon| 
&\leq \left(\int_\Omega |x|^\beta v^k |\nabla v|^2\right)^{1/2} \times \eqref{idB}^{1/2} + \eqref{idA}^{1/2} \times
\eqref{idB}^{1/2} \\& \qquad+ C_{k,\varepsilon} \left\{
\eqref{idE}^{1/2} \times\eqref{idB}^{1/2}+\eqref{idD}^{1/2} \times\eqref{idB}^{1/2}
\right\}.
\end{align*}

\Proof[ of \eqref{identity1}]
Let us now prove the second statement of Lemma~\ref{lemma4}.
One denotes by $u\in H^1_0$ a weak limit of a minimizing sequence $(u_j)_{j\in\mathbb{N}}$ for  $S_\lambda(\beta,k)$.
According to Fatou's lemma, one has $\int_\Omega |x|^\beta |u|^k |\nabla u|^2<\infty$ so
one can apply the first part of our Lemma~\ref{lemma4},  and thus
one only needs to prove the upper bound on $E_\lambda(u)$.
As $u_j$ is a minimizing sequence and $u_j\to u$ in $L^2$:
\begin{equation}\label{eq13}
S_\lambda(\beta,k) +\lambda \int_\Omega u^2 =\int_\Omega (\alpha +|x|^\beta |u_j|^k)|\nabla u_j|^2+o(1).
\end{equation}
On the other hand, the classical Brezis-Lieb Lemma $\|u\|_q^q + \|u_j-u\|_q^q=1+o(1)$ can be rewritten
\[
\|u_j-u\|_q^2 = (1-t^q)^{2/q}+o(1)
\]
and therefore, by definition \eqref{defS} of $S$:
\[
\int_\Omega |\nabla(u_j-u)|^2 \geq S \|u_j-u\|_q^2 = S(1-t^q)^{2/q}+o(1).
\]
As $u_j \rightharpoonup u$ weakly in $H^1(\Omega)$, the left-hand side expands to
\[
\int_\Omega |\nabla(u_j-u)|^2 = \int_\Omega |\nabla u_j|^2 - \int_\Omega |\nabla u|^2 + o(1)
\]
and thus one gets
\begin{equation}\label{eq14}
\int_\Omega |\nabla u_j|^2 \geq S(1-t^q)^{2/q} + \int_\Omega |\nabla u |^2 + o(1).
\end{equation}
Combining \eqref{eq13} and \eqref{eq14} gives:
\[
S_\lambda(\beta,k) +\lambda \int_\Omega u^2 \geq
\int_\Omega |x|^\beta |u_j|^k |\nabla u_j|^2 + \alpha S(1-t^q)^{2/q}+\alpha  \int_\Omega |\nabla u|^2+o(1).
\]
Fatou's lemma provides a lower bound:
\[
\int_\Omega |x|^\beta|u|^k|\nabla u|^2 \leq \operatorname{liminf} \int_\Omega |x|^\beta |u_j|^k |\nabla u_j|^2
\]
so one gets:
\[
S_\lambda(\beta,k) +\lambda \int_\Omega u^2 \geq
\int_\Omega |x|^\beta |u|^k |\nabla u|^2 + \alpha S(1-t^q)^{2/q}+\alpha  \int_\Omega |\nabla u|^2
\]
i.e. $S_\lambda(\beta,k) \geq E_\lambda(u) +\alpha  S(1-t^q)^{2/q}$.
This concludes the proof of Lemma~\ref{lemma4}.
\cqfd

\subsection{Second calculus of variations around the weak limit}\label{parcv2}

In this section, let us now explore $H^1_0(\Omega)$ around the weak limit $u$ using $v=(1+\theta)u$
within the lemma proved in the previous section.
\begin{lemma}\label{lemma5}
If $t=\|u\|_q<1$, the weak limit $u$ of the minimizing sequence $(u_j)_{j\in\mathbb{N}}$ satisfies:
\begin{equation}\label{euler-lagrange}
E_\lambda(u)+\frac{k}{2}\int_\Omega |x|^\beta |u|^k |\nabla u|^2 = \alpha S(1-t^q)^{2/q-1}\cdot t^q.
\end{equation}
\end{lemma}
This lemma is the Euler-Lagrange substitute within our reductio ad absurdum from the (not yet proven to be
bogus) assumption that $t<1$.

\Proof
The key point is that, if $t<1$, then for $\theta$ small enough one has $(1+\theta)t<1$ and thus
according to Lemma~\ref{lemma4}:
\[
S_\lambda(\beta,k)\leq E_\lambda((1+\theta)u)+\alpha S\left\{1-(1+\theta)^q t^q\right\}^{2/q}
\]
with equality when $\theta=0$.
One can thus claim that:
\[
\left.\frac{d}{d\theta} E_\lambda((1+\theta)u)+\alpha S\left\{1-(1+\theta)^qt^q\right\}^{2/q} \right|_{\theta=0}=0
\]
and a straightforward computation of the Taylor expansion of this expression with respect to $\theta$ gives:
\[
E_\lambda(u) + \frac{k}{2}\int_\Omega |x|^\beta |u|^k |\nabla u|^2 = \alpha S(1-t^q)^{2/q-1}\cdot t^q,
\]
which is exactly \eqref{euler-lagrange}.
\cqfd

\subsection{Proof of Theorem~\ref{mainThm}}\label{proofMainThm}

We are now ready to conclude the reduction ad-absurdum in which
one assumed that $t=\|u\|_q<1$.
The previous lemmas \ref{aprioriestimate}, \ref{lemma4} and \ref{lemma5} ensure respectively that:
\begin{gather}
\label{comparisonS} S_\lambda(\beta,k) < \alpha S\\
\label{exactS} S_\lambda(\beta,k)-\alpha S(1-t^q)^{2/q} = E_\lambda(u)\\
\label{euler-lagrange2}
E_\lambda(u)+\frac{k}{2}\int_\Omega |x|^\beta |u|^k |\nabla u|^2 = \alpha S(1-t^q)^{2/q-1} t^q.
\end{gather}

By definition \eqref{lambda1} of $\lambda_1(\Omega)$, one has
$\int_\Omega |\nabla u|^2\geq \lambda_1(\Omega)\int_\Omega u^2$
and thus
\begin{equation}
\label{comparisonE}
 \int_\Omega |x|^\beta |u|^k|\nabla u|^2 \leq E_\lambda(u)
 + (\lambda-\alpha \lambda_1(\Omega)) \int_{\Omega} |u|^2.
\end{equation}
Let us recall that one assumes $\lambda\leq\alpha  \lambda_1(\Omega)$ so \eqref{comparisonE}
guaranties that $E_\lambda(u)\geq0$. One can also discard the last negative term in~\eqref{comparisonE}. Combined with \eqref{exactS}, it leads us to $S_\lambda(\beta,k)>0$
as long as one assumes that $t<1$.

\medskip
Let us now combine \eqref{comparisonS}-\eqref{comparisonE} in one inequality with $S_\lambda(\beta,k)$ on both sides. More precisely, let us start by combining \eqref{euler-lagrange2} and \eqref{comparisonE}, which gives:
\[
\left(1+\frac{k}{2}\right) E_\lambda (u) \geq\alpha  S(1-t^q)^{2/q-1}t^q.
\]
Next, one substitutes the exact value for $E_\lambda(u)$ given by \eqref{exactS}:
\[
\left(1+\frac{k}{2}\right) \left(S_\lambda(\beta,k)-\alpha S(1-t^q)^{2/q}\right) \geq \alpha S(1-t^q)^{2/q-1}t^q.
\]
Then one uses \eqref{comparisonS} on both sides and $0<t<1$:
\[
\left(1+\frac{k}{2}\right) \left(S_\lambda(\beta,k)-S_\lambda(\beta,k)(1-t^q)^{2/q}\right) > S_\lambda(\beta,k)\cdot
(1-t^q)^{2/q-1}t^q.
\]
As  $S_\lambda(\beta,k)>0$,  one can simplify by $S_\lambda(\beta,k)$ and get an equivalent statement:
\[
1+\frac{k}{2} > \frac{ (1-t^q)^{2/q-1}t^q}{1-(1-t^q)^{2/q}}\cdotp
\]
We now claim that the right-hand side is an increasing function of $t$ on $(0,1)$ which is bounded from below by $q/2$ which is the limit at the origin. Therefore, one gets a contradiction as soon as:
\begin{equation}
0\leq k\leq q-2.
\end{equation}
To back-up our  claim, let us compute the derivative:
\[
\frac{d}{dt}\left(\frac{ (1-t^q)^{2/q-1}t^q}{1-(1-t^q)^{2/q}}\right) = \frac{t^{q-1}(1-t^q)^{-2+\frac{2}{q}}\left\{
q(1-(1-t^q)^{2/q})-2t^q
\right\}  }{(1-(1-t^q)^{2/q})^2}
\]
The function $q(1-(1-t^q)^{2/q})-2t^q$ itself is an increasing function of $t$ on $(0,1)$ that vanishes at $t=0$. It is
therefore positive which in turn implies that $\frac{ (1-t^q)^{2/q-1}t^q}{1-(1-t^q)^{2/q}}$ is increasing too.
The limit
\[
\lim_{t\to0^+} \frac{ (1-t^q)^{2/q-1}t^q}{1-(1-t^q)^{2/q}} = 
\lim_{x\to0^+} \frac{ x}{1-(1-x)^{2/q}} =\frac{q}{2}
\]
is therefore a lower bound for the function on $(0,1)$, which settles the claim.

\bigskip
The general conclusion of this reduction ad-absurdum that we started p.\pageref{beginContradiction}, is that $t=\|u\|_q=1$ and therefore $u$ itself is a minimizer.
This concludes the proof of Theorem~\ref{mainThm}.

\bigskip
\paragraph{Remarks on the case $\lambda>\lambda_1(\Omega)$.}
Let us briefly discuss what results survive when $\lambda$ exceeds its critical value.
One can still use~\eqref{comparisonE} and the previous inequality~:
\[
\forall t\in (0,1),\qquad  (1-t^q)^{2/q-1} t^q > \frac{q}{2}\left(1-(1-t^q)^{2/q}\right).
\]
Assuming that $S_\lambda(\beta,k)>0$ the previous computation boils down to:
\[
\left(\frac{q-2}{k}-1\right) \frac{1}{|\Omega|^{1-2/q}} \frac{S_\lambda(\beta,k)}{\lambda-\alpha \lambda_1(\Omega)}
< \frac{t^2}{1-(1-t^q)^{2/q}}\cdotp
\]
But as the right-hand side is a decreasing bijective map from $(0,1)$ to $(1,+\infty)$, one gets a restriction on $t$ of the form $t<T(\Omega,\lambda)$, as long as the left-hand side exceeds 1.
In other words, one gets the following partial result that could at least be interesting for the numerical
analysis of this problem.

\begin{prop}
Let us assume that $n\geq4$, $q=\frac{2n}{n-2}$, $0\leq k\leq q-2$, $\beta > \frac{kn}{q}+2$ and that
 \begin{equation}
\alpha  \lambda_1(\Omega)<\lambda<\alpha \lambda_1(\Omega)+\left(\frac{q-2}{k}-1\right) \frac{S_\lambda(\beta,k)}{|\Omega|^{1-2/q}}
 \qquad\text{with}\qquad S_\lambda(\beta,k)>0.
\end{equation}
Then for any minimizing sequence $u_j \in H^1_0(\Omega)$ i.e. such that
\[
\|u_j\|_{L^q}=1 \quad\text{and}\qquad E_\lambda(u_j)=S_\lambda(\beta,k)+o(1)
\]
that converges weakly to some $u\in H^1_0(\Omega)$, then either $\|u\|_{L^q}=1$ and $u$ is a minimizer,
or one has the alternative $\|u\|_{L^q}\in(0, T(\Omega,\lambda))$ where $T(\Omega,\lambda)$ is the unique
solution in $(0,1)$ of
\[
\left(\frac{q-2}{k}-1\right) \frac{1}{|\Omega|^{1-2/q}} \frac{S_\lambda(\beta,k)}{\lambda-\alpha \lambda_1(\Omega)}
= \frac{T(\Omega,\lambda)^2}{1-(1-T(\Omega,\lambda)^q)^{2/q}}\cdotp
\]
\end{prop}

\subsection{Proof of Proposition~\ref{corollary}}\label{proofMainThm2}

The previous proof holds for any minimizing sequence $u_j\in H^1_0(\Omega)$ that converges weakly to some function $u\in H^1_0(\Omega)$.
From what precedes, one can claim that $\|u\|_q=1$ and that $u$ is a minimizer.

\subsubsection{Strong convergence in $H^1(\Omega)$}

To prove Proposition~\ref{corollary}, one first needs to show that the convergence happens in the strong topology of $H^1(\Omega)$.
As $u$ is a minimizer, one has:
\begin{equation}
\int_\Omega (\alpha +|x|^\beta |u|^k)|\nabla u|^2 = \int_\Omega (\alpha +|x|^\beta |u_j|^k)|\nabla u_j|^2 + o(1)
\end{equation}
and according to Fatou's lemma
\[
\int_\Omega |\nabla u|^2 \leq \lim_{j\to+\infty} \int_\Omega |\nabla u_j|^2.
\]
If this inequality were strict, then one should also have
\[
\int_\Omega |x|^\beta |u|^k |\nabla u|^2 >   \lim_{j\to+\infty}\int_\Omega |x|^\beta |u_j|^k |\nabla u_j|^2
\]
in order for the sum of the left-hand sides to equal the sum of the right-hand sides. However, this last one violates Fatou's lemma. One must thus have,
as $\alpha>0$:
\begin{equation}
\int_\Omega |\nabla u|^2 = \lim_{j\to+\infty} \int_\Omega |\nabla u_j|^2
\end{equation}
and the strong convergence follows from the following classical trick in Hilbert spaces:
\[
\lim_{j\to+\infty}\|u_j-u\|_{H^1}^2 =  \|u\|_{H^1}^2  + \lim_{j\to+\infty}\|u_j\|_{H^1}^2 - 2\lim_{j\to+\infty} (u_j \vert u)_{H^1} = 0.
\]

\subsubsection{Euler-Lagrange equation}
The last point of Proposition~\ref{corollary} is that $u$ is a non-trivial solution of the Euler-Lagrange equation.
\begin{equation}
\begin{cases}
-\div\left(\left(\alpha +|x|^\beta |u|^k\right)\nabla u\right) +  \frac{k}{2} |x|^\beta |u|^{k-2}u |\nabla u|^2 =\lambda u + \Theta
|u|^{q-2}u\\
u_{\vert \partial\Omega}=0
\end{cases}
\end{equation}
for some $\Theta>0$. The non-trivial part is simply that $\|u\|_q=1$.

%with $\displaystyle \Theta= \frac{1}{q}\left\{\int_\Omega|\nabla u|^2+\left(\frac{k}{2}+1\right)\int_\Omega |x|^\beta |u|^k |\nabla u|^2 -\lambda \int_\Omega u^2 \right\}$. 

\bigskip
For any $\varphi \in H^1_0$ and $\theta\in\R$, the
function $\frac{u+\theta\varphi}{\|u+\theta\varphi\|_q}$ is an acceptable test function so
\[
S_\lambda(\beta,k) \leq E_\lambda\left(\frac{u+\theta\varphi}{\|u+\theta\varphi\|_q}\right).
\]
Formally, the right-hand side can be developed as a power series in $\theta$, which takes the form:
\[
E_\lambda\left(\frac{u+\theta\varphi}{\|u+\theta\varphi\|_q}\right)
= E_\lambda\left(\frac{u}{\|u\|_q}\right) + \theta \int_\Omega F_\lambda(u)\varphi + o(\theta).
\]
If $u$ is a minimizer constructed in Theorem~\ref{mainThm}, then $\|u\|_q=1$ and $E_\lambda\left(\frac{u}{\|u\|_q}\right)= S_\lambda(\beta,k)$. To satisfy the variational inequality for any small $\theta\in\R$, one must then have $F_\lambda(u)=0$, which is the Euler-Lagrange equation.
However, as $E_\lambda(v)=+\infty$ for some $v\in H^1_0(\Omega)$, one must be careful and check that the left-hand side is
indeed a $C^1$ function of $\theta$ near the origin.

\medskip
Let us restrict ourselves
to smooth, compactly supported test functions, i.e. $\varphi\in C^\infty_c(\Omega)$ ; $u$ is the minimizer constructed in Theorem~\ref{mainThm}.
In this case, as $H^1_0(\Omega)\subset L^q(\Omega) \subset L^k(\Omega)$ for $k\leq q$, one has:
\begin{align*}
\int_\Omega |x|^\beta |u+\theta \varphi|^k |\nabla u+\theta\nabla\varphi|^2
&\leq C_k \int_\Omega |x|^\beta (|u|^k + C)  (|\nabla u| ^2 +C)
\\&\leq C_k \int_\Omega |x|^\beta |u|^k |\nabla u|^2 + C_k'\|u\|_{H^1}^2 + C_k''\|u\|_{H^1}^k + C_k'''.
\end{align*}
This ensures that $E_\lambda\left(\frac{u+\theta\varphi}{\|u+\theta\varphi\|_q}\right)<\infty$ and gives
meaning to the previous formal argument.
To get the equation, the computation goes as follows:
\[
\|u+\theta\varphi\|_q^\sigma = \left(\int_\Omega |u|^q+\theta \int_\Omega |u|^{q-2} u \varphi + o(\theta) \right)^{\sigma/q}
= 1 +  \frac{\sigma \theta}{q}\int_\Omega |u|^{q-2} u \varphi  + o(\theta)
\]
and
\begin{align*}
E_\lambda\left(\frac{u+\theta\varphi}{\|u+\theta\varphi\|_q}\right) &=
\|u+\theta\varphi\|_q^{-2}\left\{
\int_\Omega \left(\alpha +\frac{|x|^\beta |u+\theta\varphi|^k}{\|u+\theta\varphi\|_q^k}\right) |\nabla u+\theta\nabla\varphi|^2 - \lambda \int_\Omega (u+\theta\varphi)^2 \right\}\\
&=  \left( 1- \frac{2\theta}{q}\int_\Omega |u|^{q-2} u \varphi \right)\left\{
\int_\Omega \alpha |\nabla u|^2 - \lambda u^2 + 2 \theta \cdot \left(\alpha  \int_\Omega \nabla u\nabla \varphi -\lambda \int_\Omega u \varphi\right)
\right.\\
& \qquad \left. + \left(1- \frac{k \theta}{q}\int_\Omega |u|^{q-2}u\varphi \right)\cdot
\left( \int_\Omega |x|^\beta (|u|^k+k\theta |u|^{k-2}u\varphi) (|\nabla u|^2+2\theta\nabla u\nabla\varphi) \right) \right\}\\
& \qquad + o(\theta).
\end{align*}
The term of order one in $\theta$ is:
\begin{align*}
&
%\left.\frac{d}{d\theta}E_\lambda\left(\frac{u+\theta\varphi}{\|u+\theta\varphi\|_q}\right)\right|_{\theta=0} &=
 2 \left(\alpha  \int_\Omega \nabla u\nabla \varphi -\lambda \int_\Omega u \varphi\right)\\
& \qquad +  2\int_\Omega |x|^\beta |u|^k \nabla u \nabla \varphi 
+k \int_\Omega |x|^\beta |u|^{k-2}u |\nabla u|^2 \varphi \\
&\qquad - \int_\Omega |u|^{q-2}u \varphi \left\{\frac{2}{q}\left(\int_\Omega \alpha |\nabla u|^2-\lambda u^2\right) +\left(\frac{k+2}{q}\right)\int_\Omega |x|^\beta |u|^k |\nabla u|^2\right\}
=\int_\Omega F_\lambda(u)\varphi.
\end{align*}
with $\frac{1}{2}F_\lambda(u)=-\alpha \Delta u -\lambda u - \div(|x|^\beta |u|^k\nabla u)+ \frac{k}{2} |x|^\beta|u|^{k-2}u|\nabla u|^2 - \Theta |u|^{q-2}u$
and
\[
\Theta=\frac{2}{q}\left(\int_\Omega \alpha |\nabla u|^2-\lambda u^2\right) +\left(\frac{k+2}{q}\right)\int_\Omega |x|^\beta |u|^k |\nabla u|^2
\]
thus $\Theta>0$ for $\lambda\leq \alpha \lambda_1(\Omega)$ and $u\not\equiv0$.
 \cqfd

\section{Generalisations and open problems}\label{generalisation}

Theorem~\ref{mainThm} remains valid for the following minimisation problem, which is more general.
\begin{thm}\label{mainThm2}
Let us consider:
\begin{equation}
S_\lambda(a)=\inf_{\substack{u\in H^1_0(\Omega)\\\norme[L^q]{u}=1}} \left\{ \int_\Omega a(x,u(x)) |\nabla u(x)|^2 dx - \lambda \int_\Omega u^2 \right\}
\end{equation}
where $a(x,s)=b_1(x)+b_2(x)|s|^k$. The functions $b_1$, $b_2$ satisfy the following assumptions:
\begin{enumerate}
\item $b_1$ has a global minimum $\alpha=b_1(x_0)$ at some point $x_0\in \Omega$, of order $\gamma>2$, i.e.:
\begin{equation}
\begin{cases}
b_1(x)=\alpha + O(|x-x_0|^\gamma),\qquad \gamma>2\\
b_1(x)\geq\alpha \quad\text{if}\quad  x\neq x_0\end{cases}
\end{equation}
\item$b_2$ is positive and has a unique zero in $\Omega$ at the same point $x_0$, of order $\beta$
\begin{equation}
\begin{cases}
b_2(x)=|x-x_0|^\beta + o(|x-x_0|^\beta)\\
b_2(x)>0 \quad\text{if}\quad x\neq x_0.
\end{cases}
\end{equation}
\end{enumerate}
One assumes restrictions on the parameters that are similar to the ones in Theorem~\ref{mainThm}:
 \begin{equation}
 0< \lambda\leq \alpha\lambda_1(\Omega), \qquad
0\leq k\leq q-2 \qquad\text{and}\qquad \beta > \frac{kn}{q}+2.
\end{equation}
Then there exists $u\in H^1_0(\Omega)$ with $\|u\|_{L^q}=1$ such that
\begin{equation}
S_\lambda(a) = \int_\Omega a(x,u(x)) |\nabla u(x)|^2 dx - \lambda \int_\Omega u^2
\end{equation}
\end{thm}
\paragraph{Remark.}
Actually, the result (and our proof) would remain valid for a general function $a(x,s)$ such that:
\begin{equation}\label{bestAssumption}
b_1^-(x)+b_2^-(x)|s|^k \leq a(x,s) \leq b_1^+(x)+b_2^+(x)|s|^k
\end{equation}
with $b_1^\pm$ and $b_2^\pm$ that satisfy similar assumptions.

\paragraph{Sketch of proof.}
Let us focus briefly on how one would adapt the previous proof to deal with this case.
The first crucial step is the comparison between $S_\lambda(a)$ and $\alpha S$ (Lemma \ref{aprioriestimate}). The remainder of $b_1$ produces an additional
term
\[
\int_\Omega |x-x_0|^\gamma \left|\nabla\left(\frac{\omega_\varepsilon}{\|\omega_\varepsilon\|_q}\right)\right|^2 = O(\varepsilon^{\gamma/2})
\]
which is negligible in comparison to $-\lambda \int_\Omega \left(\frac{\omega_\varepsilon}{\|\omega_\varepsilon\|_q}\right)^2\sim - C\lambda \varepsilon$
provided $\gamma>2$.
The second crucial step is to prove \eqref{identity1}. Provided $b_1(x)\geq b_1(x_0)=\alpha$, i.e. that the minimum at $x_0$ is a global one, one still gets :
\[
\int_\Omega b_1(x) |\nabla u_j|^2 \geq \alpha S(1-t^q)^{2/q} + \int_\Omega b_1(x) |\nabla u |^2 + o(1)
\]
instead of \eqref{eq14} and the rest of the proof remains unchanged. \cqfd

\bigskip
Further generalisations beyond Theorem~\ref{mainThm2}, or at best~\eqref{bestAssumption}, seem for now out of reach.
\begin{itemize}
\item
If $b_1$ admits a minimum of order $\gamma\leq 2$ at $x_0$, the comparison between $S_\lambda(a)$ and $\alpha S$ is not clear anymore.
\item 
If $b_1$ admits only a local minimum at $x_0$ which is not global, most of the comparisons that we used cease to be true.
\item
If the zero of $b_2$ occurs at a point $x_1$ distinct of the point $x_0$ where $b_1$ reaches its minimum, then minimizing sequences can either concentrate around
$x_0$ and one expect a behaviour similar to the model case with $\beta=0$ (i.e. the solution is in the linear regime studied in~\cite{HV} and a minimizer should exist),
or the the minimizing sequence can concentrate around $x_1$ and in that case again, it is not clear how to compare $S_\lambda(a)$ and $(\min a)\cdot S$ anymore.
\end{itemize}

\end{document}